\documentclass{article}
\usepackage{latexsym,amsmath,euscript,amsthm,amsfonts,epsfig,psfrag }
\begin{document}

\newtheorem{lemma}{Lemma}
\newtheorem{theorem}{Theorem}
\newtheorem{cor}{Corollary}
\renewcommand{\proofname}{Proof}

\newcommand{\Q}{\mathbb{Q}}
\newcommand{\R}{\mathbb{R}}
\newcommand{\C}{\mathbb{C}}
\newcommand{\N}{\mathbb{N}}
\newcommand{\Z}{\mathbb{Z}}
\def \E{{I\!\!E}}
\renewcommand{\P}{\EuScript P}
\newcommand{\Vol}{V\hspace{-2.5pt}ol}

\begin{center}
{\LARGE  \bf
 On the volume of a six-dimensional\\
\vspace{5pt} 
polytope }

\vspace{15pt} 
{\large A.~Felikson, P.~Tumarkin}

\end{center}

\vspace{9pt}

\begin{center}
\parbox{10cm}{ \small
{\bf Abstract.}
This note is a comment to the paper~\cite{MN}. That paper 
concerns with the projective surface $S$ in $\mathbb{P}^{3}$ defined by the equation
$x_{1}x_{2}x_{3}=x_{4}^{3}$. It is shown there that the evaluation of the leading term
of  the asymptotic formula 
for the number of rational points of bounded height in $S(\Q)$
is equivalent to the evaluation of the volume of some 6-dimensional polytope $\P$.
The volume of $\P$ is known from several papers; 
we calculate this volume by elementary method using symmetry of the polytope $\P$.
We also discuss a combinatorial structure of $\P$.\\

}
\end{center}

\section{Introduction}

Consider the projective surface $S$ in $\mathbb{P}^{3}$ defined by the equation
$$x_{1}x_{2}x_{3}=x_{4}^{3}.$$ 
A few years ago several authors, \cite{CD}, \cite{GH}, \cite{IK}, \cite{MN}, \cite{OP}, 
obtained an asymptotic formula for the number of rational points of bounded height in $S(\Q)$, 
as a simple illustration of a rather general conjecture discussed in \cite{AB}, \cite{CD}, \cite{OP}. 
To explicitly evaluate the leading term of the asymptotic formula obtained in \cite{CD}, \cite{MN}, \cite{OP} 
one must calculate the volume $\Vol(\P)$ of the polytope
$\P\in \R^{6}$ determined by the inequalities

$$x_{12}+x_{13}+2(x_{21}+x_{31})\le 1,$$ 
$$x_{21}+x_{23}+2(x_{12}+x_{32})\le 1,$$  
$$x_{31}+x_{32}+2(x_{13}+x_{23})\le 1,$$   
$$x_{ij}\ge 0,$$
where 
$$(x_{12},x_{23},x_{31},x_{13},x_{32},x_{21} )$$
are Cartesian coordinates in $\R^6$.

As it is shown in papers \cite{CD}, \cite{MN} and \cite{OP},  
\begin{equation}
\Vol(\P)=\frac{1}{4\cdot 6!}=\frac{1}{2880}.
\end{equation}

The goal of this note is to prove formula (1) by a direct elementary method, using symmetry of the polytope $\P$. 
We dissect $\P$ into $6$ congruent parts and decompose one of these parts into $6$ simplices. Then we show that 
the volumes of these simplices sum up to $\frac{1}{24}\Vol(\Delta)$, where
$\Delta$ is the standard six-dimensional simplex determined by the inequalities
$$
x_{12}+x_{23}+x_{31}+x_{13}+x_{32}+x_{21}\le 1,\quad x_{ij}\ge 0,\qquad  i\ne j,\quad 1\le i,j\le 3.
$$


We use the notation and terminology of \cite{Grunbaum}. 
In particular,
given a subset $X$ of $\R^{n},\;n\in \N$, we write $conv\; X$ and $\partial X$ for the convex hull and the boundary of $X$ respectively. 

\vspace{8pt}
{\bf Acknowledgment.}
We would like to thank  B.~Z.~Moroz  who brought the problem to our attention,
initiated  numerous stimulating discussions and convinced us to write down this note.

\section{Volume of a rational polytope}

Let $\E^n$ be the $n$-dimensional  Euclidean space with Cartesian coordinates. 
Let $\bf v_0,\bf v_1,\dots,\bf v_k$ be rational points in $\E^n$, i.e. the coordinates $(v_{i1},\dots,v_{in})$ of $\bf v_i$ are rational numbers for all $i=0,\dots, k$.   
Suppose that the points  $\bf v_0,\bf v_1,\dots,\bf v_k$ are in the convex position, that is no of these points belongs to the
convex hull of the others. 
Suppose that $k\ge n$. Then the  set  $conv\{\bf v_0,\dots,\bf v_k\}$ is an $n$-dimensional polytope with non-zero  volume $\Vol(\bf v_0,\dots,\bf v_k)$.
In this section we describe how to find $\Vol(\bf v_1,\dots,\bf v_k)$.


Denote by $det(\bf v_1,\dots,{\bf  v_n})$ the determinant 
$$
det(\bf  v_1,\dots,\bf  v_n)=
\begin{vmatrix}
v_{11}  &v_{12} & \ldots & v_{1n} \\
v_{21}  &v_{22} & \ldots & v_{2n} \\
\vdots  &\vdots & \ddots & \vdots \\
v_{n1}  &v_{n2} & \ldots & v_{nn} \\
\end{vmatrix}.
$$

If $k=n$, the polytope  $conv\{\bf v_0,\bf v_1,\dots,\bf v_k\}$ is a simplex (may be degenerate). 
Recall that 
$$
\Vol({\bf v_0},{\bf v_1},\dots,{\bf v_n})=
\frac{1}{n!}\cdot 
det({\bf v_1}-{\bf v_0},{\bf v_2}- {\bf v_0},\dots,{\bf  v_n}- {\bf v_0})
$$
(see for example~\cite{Berger}, Chapter 9).
%
\vspace{6pt}

\noindent
{\bf Definition.}
{\it A triangulation} of a polytope is a decomposition of this polytope into finite number of simplices $T_1,\dots,T_r$, such that
if $T_i\bigcap T_j \ne \emptyset$, then   $T_i\bigcap T_j$ is a face of both $T_i$ and $T_j$.

\vspace{6pt}

The formula above reduces the problem to the question of triangulation of the polytope $conv\{\bf v_0,\bf v_1\dots,\bf v_k\}$.

Lemmas~\ref{half} through~\ref{vis} are either well-known or elementary.


\begin{lemma}
\label{half}
Let $\bf v_i$, $1\le i\le n$ be $n$ points in general position in $\E^n$
(i.e. no $k$ of these points belong to $(k-2)$-plane).
Let $H$ be the hyperplane spanned by $\bf v_1,\dots,\bf v_n$. 
Let $\bf a$ and $\bf b$ be two points in $\E^n$.
The points $\bf a$ and $\bf b$ belong to different open half-spaces with respect to $H$
if and only if 
%
%
$$det({\bf v_1}-{\bf a},{\bf v_2}-{\bf a},\dots,{\bf v_{n}}- {\bf a})\cdot det({\bf v_1}-{\bf b},{\bf v_2}-{\bf b},\dots,{\bf v_{n}}- {\bf b})<0.$$ 

\end{lemma}

%

\begin{lemma}
\label{2M}
Let 
$$S:=\{Q, Q_{i}\mid\;1\le i\le m\}$$ 
be a set of $m+1$ polytopes satisfying the following conditions:

a) $Q_{i}\subseteq Q$ for $1\le i\le m$;

b) each facet of $Q_{i}$ is a facet of one the polytopes belonging to the set

$Q\setminus \{Q_{i}\}$ for $1\le i\le m$.
 
Then
$$Q=\cup_{i=1}^{m} Q_{i}$$ 
\end{lemma}

%


\noindent
{\bf Definition.} 
Let $P$ be a polytope in $\E^n$.
A point ${\bf y}\in P$ is  {\it visual} from the point ${\bf x}\notin P$, if
$[{\bf x},{\bf y}]\cap P={\bf y}$ (where $[{\bf x},{\bf y}]$ is the segment 
$\{\lambda {\bf a} + (1 - \lambda){\bf b}\mid\;0\le \lambda\le 1\}$. 
A facet $f$ of $P$ is {\it visual} from the point ${\bf x}\notin P$,
if each point of $f$ is visual from ${\bf x}$. 

If $P$ is convex any facet containing at least one inner visual point is visual.

\begin{lemma}
\label{vis}
Let  $P$ be a convex polytope, and ${\bf x} \notin P$ be a point.
Any point ${\bf y} \in P$ visual from $\bf x $ is contained in at least one visual facet of $P$.

\end{lemma}

%


%
%
%

\begin{theorem}
\label{triang}
Let $T=\bigcup_{i=1}^m T_i$ be a triangulation of  $conv\{\bf v_1,\dots,\bf v_k\}$.
Let ${\bf v_0}\notin conv\{{\bf v_1},\dots,{\bf v_k}\}$.
Let $f_1,\dots,f_r$ be all facets of $T_1,\dots,T_m$ visual from $\bf v_0$.
Denote $T_{m+j}:=conv\; {\bf v_0}\cup f_j$,\quad $1\le j \le r$.

Then $T'=\bigcup_{i=1}^{m+r} T_i$ is a triangulation of  $conv(\bf v_0,\bf v_1,\dots,\bf v_k)$.

\end{theorem}

The proof is obvious.

\vspace{6pt}

Theorem~\ref{triang} suggests an {\bf algorithm} 
of triangulation of  $conv\{\bf v_0,\bf v_1,\dots,\bf v_k\}$:
Consider any $n+1$  non-coplanar points (without loss of generality we may assume that these points are $\bf v_0,\bf v_1,\dots,\bf v_{n}$).
To check that these points are non-coplanar make sure that $det({\bf v_1}-{\bf v_{0}},\dots,{\bf  v_n}-{\bf v_{0}})\ne 0$. 
Then $T_1=conv\{\bf v_0,\bf v_1,\dots,\bf v_n\}$ is a simplex and $T_1$ is triangulated by itself.
Use Theorem~\ref{triang} to successively triangulate  $conv\{\bf v_0,\bf v_1,\dots,\bf v_{i}\}$ for $i=n+1,n+2,\dots,k$.


\section{Combinatorial structure of the polytope $\P$} 

Denote by $\alpha_1$, $\alpha_2$ and $\alpha_3$ the hyperplanes 
$$x_{12}+x_{13}+2(x_{21}+x_{31})=1,$$ 
$$x_{21}+x_{23}+2(x_{12}+x_{32})=1,$$  
$$x_{31}+x_{32}+2(x_{13}+x_{23})=1$$  respectively. 
Denote by $\pi_{ij}$ the hyperplane $x_{ij}=0$ ($i,j=1,2,3$,\quad $i\ne j$).
Denote by $\sigma$ the hyperplane $x_{12}+x_{23}+x_{31}+x_{13}+x_{32}+x_{21}=1$.

It is easy to see, that if $X\in \sigma$ and $X\in \P$  then $X\in \alpha_i$ for $i=1,2,3$.
Therefore, $\sigma$ is not a facet of $\P$.
Thus, $\P$ is bounded by $\alpha_i$, $i=1,2,3$  and six facets $\pi_{ij}$.
Any vertex of $\P$ is an intersection of at least 6 facets. 
Thus, any vertex of $\P$ is a subject to at least 6 equations under consideration. 
A straightforward calculation shows that $\P$ is a convex hull of 21 points.
We list these points in Table~\ref{vertices}.
The right column of the table shows, to which of the hyperplanes 
$\alpha_{j},\;1\le j\le 3,$ if any, the corresponding point belongs.
Note, that the 
indices indicate non-zero coordinates.

\begin{table}[h!]
\begin{center}
\label{vertices}
\caption{Vertices of $\P$.}
\vspace{6pt}

\begin{tabular}{|c|c|ccc|}
\hline
Notation & 
\begin{tabular}{c}
\vphantom{$\int\limits^A$}
Coordinates \\
\vphantom{$\int\limits_A$}
$(x_{12},x_{23},x_{31},x_{13},x_{32},x_{21} )$\\
\end{tabular}
& 
$\alpha_1$ & $\alpha_2$ & $\alpha_3$ \\

\hline
\vphantom{$\int\limits_A^A$}
$O$ & $(0,0,0,0,0,0)$ & $-$ & $-$ & $-$\\
\vphantom{$\int\limits_A$}
$Q_{23}^{21}$ & $(0,\frac{1}{2},0,0,0,\frac{1}{2})$ & + & + & +\\
\vphantom{$\int\limits_A$}
$Q_{32}^{31}$ & $(0,0,\frac{1}{2},0,\frac{1}{2},0)$ & + & + & +\\
\vphantom{$\int\limits_A$}
$Q_{13}^{12}$ & $(\frac{1}{2},0,0,\frac{1}{2},0,0)$ & + & + & +\\
\vphantom{$\int\limits_A$}
$P_{12}$ &  $(\frac{1}{2},0,0,0,0,0)$ & $-$ & + & $-$\\
\vphantom{$\int\limits_A$}
$P_{23}$ &  $(0,\frac{1}{2},0,0,0,0)$ & $-$ & $-$ & +\\
\vphantom{$\int\limits_A$}
$P_{31}$ &  $(0,0,\frac{1}{2},0,0,0)$ & + & $-$ & $-$\\
\vphantom{$\int\limits_A$}
$P_{13}$ &  $(0,0,0,\frac{1}{2},0,0)$ & $-$ & $-$ & +\\
\vphantom{$\int\limits_A$}
$P_{32}$ &  $(0,0,0,0,\frac{1}{2},0)$ & $-$ & + & $-$\\
\vphantom{$\int\limits_A$}
$P_{21}$ &  $(0,0,0,0,0,\frac{1}{2})$ & + & $-$ & $-$\\
\vphantom{$\int\limits_A$}
$R_1$     &  $(\frac{1}{3},\frac{1}{3},\frac{1}{3},0,0,0)$ & + & + & +\\
\vphantom{$\int\limits_A$}
$R_2$     &  $(0,0,0,\frac{1}{3},\frac{1}{3},\frac{1}{3})$ & + & + & +\\
\vphantom{$\int\limits_A$}
$T_{12}^{21}$ &  $(\frac{1}{3},0,0,0,0,\frac{1}{3})$ & + & + & $-$\\
\vphantom{$\int\limits_A$}
$T_{23}^{32}$ &  $(0,\frac{1}{3},0,0,\frac{1}{3},0)$ & $-$ & + & +\\
\vphantom{$\int\limits_A$}
$T_{31}^{13}$ &  $(0,0,\frac{1}{3},\frac{1}{3},0,0)$ & + & $-$ & +\\
\vphantom{$\int\limits_A$}
$V_{21}^{32}$ &  $(0,0,0,0,\frac{1}{4},\frac{1}{2})$ & + & + & $-$\\
\vphantom{$\int\limits_A$}
$V_{12}^{31}$ &  $(\frac{1}{2},0,\frac{1}{4},0,0,0)$ & + & + & $-$\\
\vphantom{$\int\limits_A$}
$V_{13}^{21}$ &  $(0,0,0,\frac{1}{2},0,\frac{1}{4})$ & + & $-$ & +\\
\vphantom{$\int\limits_A$}
$V_{31}^{23}$ &  $(0,\frac{1}{4},\frac{1}{2},0,0,0)$ & + & $-$ & +\\
\vphantom{$\int\limits_A$}
$V_{32}^{13}$ &  $(0,0,0,\frac{1}{4},\frac{1}{2},0)$ & $-$ & + & +\\
\vphantom{$\int\limits_A$}
$V_{23}^{12}$ &  $(\frac{1}{4},\frac{1}{2},0,0,0,0)$ & $-$ & + & +\\

\hline
\end{tabular}

\end{center}
\end{table}

\begin{figure}[!h]
\begin{center}
\psfrag{a}{$\alpha_1$}
\psfrag{b}{$\alpha_2$}
\psfrag{g}{$\alpha_3$}
\psfrag{12}{$\pi_{12}$}
\psfrag{23}{$\pi_{23}$}
\psfrag{31}{$\pi_{31}$}
\psfrag{13}{$\pi_{13}$}
\psfrag{32}{$\pi_{32}$}
\psfrag{21}{$\pi_{21}$}
\epsfig{file=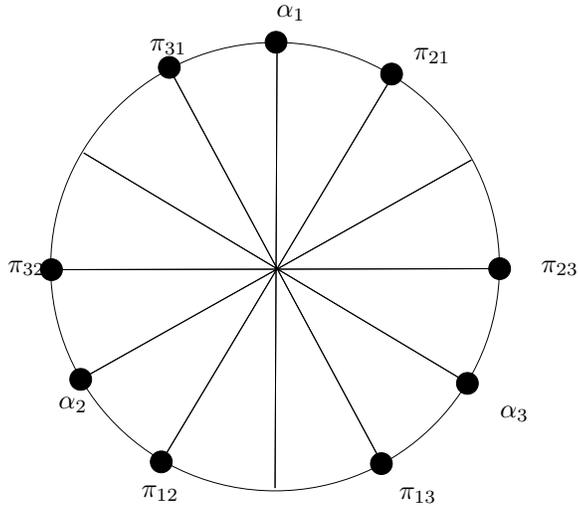,width=0.6\linewidth}
\end{center}
\caption{Gale diagram of $\P$.}
\label{gale}
\end{figure}

We present below the {\it Gale diagram} of $\P$, describing the combinatorial type of the polytope 
(see Figure~\ref{gale}).

$\P$ is a 6-dimensional polytope with 9 facets. 
The combinatorics of a convex $n$-polytope with $n+3$ facets can be described by 
2-dimensional {\it Gale diagram} (see~\cite{Grunbaum}).
 This consists of $n+3$ points  
$a_1,\dots,a_{n+3}$ of  unit circle in $\R^2$
centered  at the origin.
The combinatorial type of a convex polytope can be read off from 
the Gale diagram in the following way.
Each point $a_i$ corresponds to the facet 
$f_i$ of $\P$. For any subset $J$ of the set of
facets of $\P$ the intersection of facets 
$\{f_j | j\in J \}$ 
is a face of $\P$ if and only if the origin is contained in the set  
$conv\{a_{j} | j\notin J\}$.


\section{Triangulation of the polytope $\P$}

In this section we use the group of symmetries of $\P$ to find a nice triangulation of $\P$ containing relatively small
number of simplices.

Let 
$$\P_{j}:= conv \{O\cup \alpha_{j}\},\;1\le j\le 3.$$

\begin{lemma}
\label{3}
$Vol(\P_1)=Vol(\P_{2})=Vol(\P_{3})=\frac{1}{3}Vol(\P).$
\end{lemma}

\begin{proof}
A permutation $\varphi$ of coordinates 

$$\varphi=(x_{12}x_{23}x_{31})(x_{13}x_{21}x_{32})$$

\noindent
induces an order 3 orthogonal transformation $I_{\varphi}: \R^{6}\rightarrow \R^{6}$,
such that $I_{\varphi}(\P_{1})=\P_{2}$,
$I_{\varphi}(\P_{2})=\P_{3}$ and $I_{\varphi}(\P_{3})=\P_{1}$. Therefore,

\begin{equation}
Vol(\P_{1})=Vol(\P_{2})=Vol(\P_{3}).
\end{equation}

By construction, the set $\{\P, \P_{i}\mid\;1\le i\le 3\}$
satisfies the conditions of Lemma~\ref{2M}, hence, 
\begin{equation}
\P_{1}\cup \P_{2}\cup \P_{3}=\P.
\end{equation}
Moreover, if $1\le i<j\le 3$, then $\P_{i}\cap \P_{j}$ is a common facet of
$\P_{i}$ and $\P_{j}$; consequently,
\begin{equation}
Vol(\P_{i}\cap \P_{j})=0,\;1\le i<j\le 3.
\end{equation}
The lemma follows from (2)--(4).

\end{proof}

In view of Lemma~\ref{3}, it suffices to calculate the volume of $\P_{1}$. By construction,
$\P_{1}$ is a convex hull of the points
$$ Q_{23}^{21}, T_{12}^{21},  P_{21}, R_1, V_{21}^{32}, V_{13}^{21},$$  
$$ Q_{32}^{31}, T_{13}^{31}, P_{31}, R_2,  V_{31}^{23}, V_{12}^{31},$$ 
$$ Q_{13}^{12}, O.$$
Note that all vertices of $\P_1$ except for $Q_{12}^{13}$ and   $O$
split into pairs of similar vertices denoted by one letter with different indices. 
This will be used in the proof of Lemma~\ref{2}.

Denote by $\theta_1$ and $\theta_2$ respectively the hyperplanes
$$ x_{12}-2x_{23}+x_{31}=x_{13}+x_{32}-2x_{21}$$
and  
$$ x_{12}+x_{23}-2x_{31}=x_{13}-2x_{32}+x_{21}.$$ 

It is easy to check that $\P_1$ has the following facets:
$\pi_{ij}, \alpha_1, \theta_1$ and $\theta_2$. 


Denote by $\delta$ the hyperplane 
$$x_{12}+x_{23}+x_{31}=x_{13}+x_{32}+x_{21}.$$

\begin{lemma}
\label{2}
$\delta$ divides  $\P_1$ into two congruent parts. 

\end{lemma}

\begin{proof}
It is easy to check that $P_{31}P_{21}$ is the only edge of $\alpha_1$ which is intersected by $\delta$ transversally.

Let $K_{21}^{31}:=\delta \cap P_{31}P_{21}=(0,0,\frac{1}{4},0,0,\frac{1}{4})$. 
Then $\delta$ divides  $\P_1$ into two polytopes $\P_{1}^1$ and $\P_{1}^2$,
where  $\P_{1}^1$ is a convex hull of 
\begin{center}
$ P_{31}, R_1, V_{31}^{23}, V_{12}^{31}$ and 
$ O,Q_{23}^{21}, Q_{32}^{31}, Q_{13}^{12}, T_{12}^{21}, T_{13}^{31}, K_{21}^{31}$
\end{center}
and $\P_{1}^2$ is a convex hull of 
\begin{center}
$P_{21}, R_2, V_{21}^{32}, V_{13}^{21}$ and 
$O,Q_{23}^{21}, Q_{32}^{31}, Q_{13}^{12}, T_{12}^{21}, T_{13}^{31}, K_{21}^{31}.$
\end{center} 
Note that  $\P_1^1$ and  $\P_1^2$ differ by 4 vertices only.

Consider an involution $I_{\psi}$ induced by the following permutation of coordinates:
$\psi=(x_{12}x_{13})(x_{31}x_{21})(x_{32}x_{23})$.
%
Clearly, $I_{\psi}$ preserves $\alpha_1$ and $\delta$. More precisely,    
$I_{\psi}$ fixes the points  $Q_{12}^{13}$ and $O$
and interchanges six pairs of other points spanning $\P_1$.
Thus,  $I_{\psi}$  interchanges $\P_{1}^1$ and $\P_{1}^2$.
Since  $I_{\psi}$  is an isometry,  $\P_{1}^1$ is congruent to $\P_{1}^2$.

\end{proof}

Lemmas~\ref{3} and~\ref{2} show that the permutation group $S_3$ acts on $\P$.
Moreover, $S_3$ is a group of symmetries of $\P$ and    $\P_{1}^1=\P/S_3$. 

Lemma~\ref{2} implies that $Vol(\P_{1}^1)=Vol(\P_{1}^2)$.
By Lemmas~\ref{3} and~\ref{2} it is sufficient to triangulate  $\P_{1}^1$.
This polytope is bounded by $\pi_{ij}$, $\alpha_1$, $\theta_1$ , $\theta_2$ and $\delta$. 

We decompose  $\P_{1}^1$ into 6 simplices $\Delta_1,...,\Delta_6$ (see Table~\ref{tri}).
We represent these simplices by their vertices.
Each facet of $\Delta_i$ belongs either  to some $\Delta_j$ ($i\ne j$) or
to one of the hyperplanes bounding  $\P_{1}^1$.
Note, that a facet $f$  of a simplex corresponds to the vertex $v(f)$ opposite to this facet and vice versa.
For each vertex $v(f)$ of  $\Delta_i$ we indicate the simplex $\Delta_j$ or the facet of  $\P_{1}^1$ 
containing the facet $f$ (see the third column).

\begin{table}
\label{tri}
\caption{Triangulation of $\P_{1}^1$. }
\vspace{6pt}
\begin{tabular}{|c|c|c|c|}

\hline
&&&\\
&Vertices & Neighbors & Volume\\
&&&\\
\hline
&&&\\
$\Delta_1$ & $O,P_{31},Q_{32}^{31},Q_{23}^{21},T_{31}^{13},T_{21}^{12},K_{21}^{31}$ &  
             $\alpha_1, \delta, \pi_{32}, \pi_{23}, \pi_{13}, \pi_{12}, \Delta_2$   & 
$\frac{1}{9}\cdot \frac{2}{4^3}\frac{1}{6!} $\\ 
&&&\\
\hline
&&&\\
$\Delta_2$ & $O,P_{31},Q_{32}^{31},Q_{23}^{21},T_{31}^{13},T_{21}^{12},Q_{12}^{13}$ &  
             $\alpha_1, \delta, \pi_{32}, \pi_{23}, \Delta_3, \Delta_4, \Delta_1$   & 
$\frac{1}{9}\cdot \frac{1}{4^2}\frac{1}{6!} $\\ 
&&&\\
\hline
&&&\\
$\Delta_3$ & $O,P_{31},Q_{32}^{31},Q_{23}^{21},R_1,T_{21}^{12},Q_{12}^{13}$ &  
             $\alpha_1,\theta_2, \pi_{32}, \Delta_5, \Delta_2, \Delta_4, \pi_{13}$   &
$\frac{1}{9}\cdot \frac{2}{4^2}\frac{1}{6!} $\\ 
&&&\\
\hline
&&&\\
$\Delta_4$ & $O,P_{31},Q_{32}^{31},Q_{23}^{21},T_{31}^{13},R_1,Q_{12}^{13}$ &  
             $\alpha_1, \theta_1, \pi_{32}, \pi_{21}, \Delta_3, \Delta_2, \Delta_6$   & 
$\frac{1}{9}\cdot \frac{1}{4^2}\frac{1}{6!} $\\ 
&&&\\
\hline
&&&\\
$\Delta_5$ & $O,P_{31},Q_{32}^{31},V_{12}^{31},R_1,T_{21}^{12},Q_{12}^{13}$ &  
             $\alpha_1, \theta_2, \pi_{32}, \Delta_3, \pi_{23}, \pi_{21},  \pi_{13}$   &
$\frac{1}{9}\cdot \frac{1}{4^2}\frac{1}{6!}$ \\ 
&&&\\
\hline
&&&\\
$\Delta_6$ & $O,P_{31},Q_{32}^{31},Q_{23}^{21},T_{31}^{13},R_1,V_{31}^{23}$ &  
             $\alpha_1, \theta_1, \pi_{32}, \pi_{21}, \pi_{13}, \pi_{12}, \Delta_4$   & 
$\frac{1}{9}\cdot \frac{2}{4^3}\frac{1}{6!}$ \\ 
&&&\\
\hline 
\end{tabular}

\end{table}

\begin{theorem}
\label{vol}
1) $\P_{1}^1$ is  triangulated by $\Delta_1,\dots,\Delta_6$.\\
2) $\Vol(\P)=\frac{1}{4}\frac{1}{6!}=\frac{1}{4} \Vol(\Delta).$

\end{theorem}

\begin{proof}
The set vertices of $\Delta_1,\dots,\Delta_6$ coincides with the set of vertices of $\P_{1}^1$.
Hence, the union of $\Delta_i$ ($i=1,...,6$) lies inside  $\P_{1}^1$.
The intersection $\Delta_i\cap \Delta_j$ is a face of $\Delta_i$ and $\Delta_j$.
It is left to show that any point of  $\P_{1}^1$ belongs to some of $\Delta_i$. 
The third column of Table~\ref{tri} shows that each facet of $\Delta_i, i=1,\dots,6$ either belongs to 
$\partial \P_{1}^1$ or to some $\partial \Delta_j$, $j\ne i$.
By Lemma~\ref{2M}, $\P_{1}^1=\cup_{i=1}^6\Delta_i$ and the first statement is proved.

By lemmas~\ref{3} and \ref{2},
$$
\Vol(\P)=6\cdot \Vol(\P_{1}^1)=
6\sum_{i=1}^6 \Vol(\Delta_i).$$
In view of triangulation shown in Table~\ref{tri},
this equals to $\frac{1}{4} \cdot \frac{1}{6!}$.
This proves the first equality. The second equality is trivial.

\end{proof}

%
%
%


\vspace{25pt}
\noindent
Independent University of Moscow, Russia \\
Max-Planck Institut f\"ur Mathematik Bonn, Germany\\
e-mail: \phantom{ow} felikson@mccme.ru\qquad pasha@mccme.ru

\end{document}